\documentclass{amsart}
\usepackage{amsmath}
\usepackage{amscd}
\usepackage{amssymb}
\usepackage{mathrsfs, xcolor}

\newcommand{\mrm}{\mathrm}
\newcommand{\cal}{\mathcal}
\renewcommand{\div}{\operatorname{div}}

\newcommand{\bC}{{\Bbb C}}

\newcommand{\bR}{{\Bbb R}}

\newcommand{\cL}{{\cal L}}

\newcommand{\RP}{\mathbb{RP}}
\newcommand{\CP}{\mathbb{CP}}

\def\be{\begin{equation}}
	\def\ee{\end{equation}}

\def\al{\aligned}
\def\eal{\endaligned}

\def\d{\nabla}

\DeclareMathOperator{\Ker}{Ker}
\newtheorem{theorem}{Theorem}[section]
\newtheorem{theorem/definition}{Theorem/Definition}[section]
\newtheorem{proposition}{Proposition}[section]
\newtheorem{lemma}{Lemma}[section]

\theoremstyle{remark}
\newtheorem{remark}{Remark}[section]

\theoremstyle{definition}

\numberwithin{equation}{section}

\begin{document}
\title
{Linear stability of compact shrinking Ricci solitons}

\author{Huai-Dong Cao and Meng Zhu}
\thanks{The first author was partially supported by a Simons Foundation Collaboration Grant. % (\#586694 HC).
The second author was
partially supported by  NSFC Grant No. 11971168, Shanghai Science and Technology Innovation Program Basic Research Project STCSM 20JC1412900, and the Science and Technology Commission of Shanghai Municipality No. 22DZ2229014.}
\address{Department of Mathematics\\Lehigh University \\Bethlehem, PA 18015, USA} \email{huc2@lehigh.edu}

\address{School of Mathematical Sciences,  Key Laboratory of MEA (Ministry of Education) \& Shanghai Key Laboratory of PMMP,  East China Normal University, Shanghai 200241, China}
\email{mzhu@math.ecnu.edu.cn}

\date{}

\begin{abstract}
In this paper, we continue investigating the second variation of Perelman's $\nu$-entropy for compact shrinking Ricci solitons. In particular, we improve some of our previous work in \cite{CZ12}, as well as the more recent work in \cite{MR}, and obtain a necessary and sufficient condition for a compact shrinking Ricci soliton to be  linearly stable. % in terms of the least eigenvalue of a certain Lichnerowicz type operator associated with the stability operator of the second variation.
Our work also extends similar results of Hamilton, Ilmanen and the first author in \cite{CHI04} (see also \cite{CH15}) for positive Einstein manifolds to the compact shrinking Ricci soliton case.
 \end{abstract}
\maketitle

\section{Introduction} \label{sec:1}
\hspace{.75cm}

This is a sequel to our previous paper \cite{CZ12}, in which we derived the second variation formula of Perelman's $\nu$-entropy for compact shrinking Ricci solitons and obtained certain necessary condition for the linear stability of compact Ricci shrinkers.

Recall that a complete Riemannian manifold $(M^n, g)$ is called a {\it shrinking Ricci soliton} if there exists a smooth vector field $V$ on $M^n$ such that the Ricci tensor $Rc$
of the metric $g$ satisfies the equation
\begin{equation*}
Rc+ \frac 12 \mathscr{L}_V g=\frac 1 {2\tau} g,
\end{equation*}
where $\tau>0$ is a constant and $\mathscr{L}_V g$ denotes the Lie derivative of $g$ in the direction of $V$. If $V$ is the gradient vector field $\d f$ of a smooth function $f$, then we have a {\it gradient shrinking Ricci soliton} given by
\begin{equation} \label{eq1.1}
Rc+\nabla^2 f=\frac 1 {2\tau} g,
\end{equation}
for some constant $\tau>0$. Here, $\nabla^2f$ denotes the Hessian of $f$, and $f$ is called a {\it
potential function} of the Ricci soliton.  Clearly, when $f$ is a constant
we have an Einstein metric of positive scalar curvature. Thus,
gradient shrinking Ricci solitons include positive Einstein manifolds as a special case.  In the following, we use $(M^n, g, f)$ to denote a gradient shrinking Ricci soliton.

Gradient shrinking Ricci solitons are self-similar solutions to Hamilton's Ricci flow, and often arise as Type I singularity models in the Ricci flow as shown by Naber \cite{Na10}, Enders-M\"uller-Topping \cite{EMT11} and Cao-Zhang \cite{CxZq11}; see also Zhang \cite{Zha22}. As such, they play a significant role in the study of the formation of singularities in the Ricci flow and its applications. Therefore, it is very important to either classify, if possible,  or understand the geometry of gradient shrinking Ricci solitons.

Hamilton \cite{Ha95F} showed that any 2-dimensional complete gradient shrinking Ricci soliton is isometric to either $\mathbb S^2$, or $\RP^2$, or the Gaussian shrinking soliton on $\bR^2$. In dimension $n=3$, by using the Hamilton-Ivey curvature pinching, Ivey \cite{Iv} proved that a compact shrinking soliton must be a spherical space form $\Bbb S^3/\Gamma$. Furthermore, for $n=3$, a complete classification follows from the works of Perelman \cite{P2}, Naber \cite{Na10}, Ni-Wallach \cite{NW08}, and Cao-Chen-Zhu \cite{CCZ08} that any three-dimensional complete gradient shrinking Ricci soliton is either isometric to the Gaussian soliton $\bR^3$ or a finite quotient of either $\mathbb S^3$ or $\mathbb S^2\times \bR$.

However, in dimension $n\ge 4$, there do exist non-Einstein and non-product gradient shrinking Ricci solitons.  Specifically,  in
dimension $n=4$,  Koiso \cite{Ko} and the first author
\cite{Cao94} independently constructed a gradient K\"ahler-Ricci
shrinking soliton on $\CP^2\#(-\CP^2)$, and Wang-Zhu
\cite{WZ} found another one on $\CP^2\#(-2\CP^2)$. In the noncompact
case, Feldman-Ilmanen-Knopf \cite{FIK} constructed a
$U(2)$-invariant gradient shrinking K\"ahler-Ricci soliton on the
tautological line bundle $\mathcal{O}(-1)$ of $\CP^{1}$, i.e., the
blow-up of $\bC^{2}$ at the origin. Very recently, a noncompact toric gradient shrinking K\"ahler-Ricci soliton on the blowup of $\CP^1\times \bC$ at one point was found by Bamler-Cifarelli-Conlon-Deruelle \cite{BCCD22}. These are the only known examples of nontrivial (i.e., non-Einstein) and non-product complete shrinking Ricci solitons in dimension $4$ so far.  We remark that the constructions in \cite{Cao94, FIK, Ko, WZ} all extend to higher dimensions. For additional examples in higher dimensions, see, e.g.,  Angenent-Knopf \cite{AK22}, Dancer-Wang \cite{DW11}, Futaki-Wang \cite{FW11}, and Yang \cite{Yb12}.

Ricci solitons can be viewed as fixed points of the Ricci flow, as
a dynamical system on the space of Riemannian metrics modulo
diffeomorphisms and scalings.
 In \cite{P1}, Perelman introduced the $\mathcal{W}$-functional
$$
\mathcal{W}(\hat{g},\hat{f}, \hat{\tau})=\int_M [\hat{\tau} (\hat{R}+|\nabla \hat{f}|^2)+\hat{f}-n](4\pi
\hat {\tau})^{-\frac{n}{2}}e^{-\hat{f}}d\hat{V},
$$
on any compact manifold $M^n$, where $\hat{g}$ is a Riemannian
metric on $M$, $\hat{R}$ is its scalar curvature, $\hat{f}$ is any smooth function on $M^n$,
and ${\hat\tau}>0$ is a positive parameter. The associated
$\nu$-entropy is defined by
$$
\nu(\hat{g})=\inf\{\mathcal W(\hat{g},\hat{f},\hat{\tau}): \hat{f}\in C^\infty(M), \hat{\tau}>0,
(4\pi\hat{\tau})^{-\frac{n}{2}}\int_M e^{-\hat{f}}d\hat{V}=1\},
$$
which is always attained by some $\hat{f}$ and $\hat{\tau}$.
Furthermore, Perelman showed that the $\nu$-entropy is monotone increasing under
the Ricci flow, and its critical points are precisely given by
gradient shrinking Ricci solitons $(M^n, g, f)$ satisfying (1.1). In particular, it follows that all compact shrinking Ricci solitons are necessarily gradient ones. %, a fact shown by Perelman \cite{P1}.

By definition, a compact shrinking Ricci soliton $(M^n, g, f)$ is {\it linearly stable} (or $\nu$-stable)  if % it is a local maximum point of the $\nu$ functional, i.e.,
the second variation of the $\nu$-entropy  is nonpositive at $g$. % The study of the linear stability of compact gradient Ricci solitons was initiated by Hamilton, Ilmanen and the first author in 2004. Furthermore, it is called {\it neutrally linearly stable} if the second variation is nonpositive but not strictly negative.
In \cite{CHI04}, Hamilton, Ilmanen and the first author initiated the study of linear stability of compact shrinking Ricci solitons. They obtained the second variation formula of Perelman's $\nu$-entropy for positive Einstein manifolds and investigated
their linear stability. Among other results, they showed that, while the round sphere
$\Bbb{S}^ n$ is linearly stable and the complex projective space $\CP^n$ is neutrally linearly stable\footnote{Recently, Knopf and Sesum \cite{KS} showed that $\CP^n$ is not a local maximum of the $\nu$-entropy, hence is dynamically unstable as first shown by Kr\"oncke \cite{Kk}.}, many known positive Einstein manifolds are unstable. % for the Ricci flow so that generic perturbations acquire higher $\nu$-entropy and thus can never return near the original metric.
In particular,  all product Einstein manifolds and Fano K\"ahler-Einstein manifolds with Hodge number $h^{1,1}>1$ are unstable.  More recently, a complete description of the linear stability (or instability) of irreducible symmetric spaces of compact type was provided by C. He and the first author \cite{CH15}.
Meanwhile, in \cite{CZ12}, we derived the second variation formula of Perelman's $\nu$-entropy for  compact shrinking Ricci solitons which we now recall.

Let $(M^n, g, f)$ be a compact shrinking Ricci soliton satisfying \eqref{eq1.1} and $\mrm{Sym}^2(T^*M)$ denote the space  of symmetric (covariant) 2-tensors on $M$. For  any $h=h_{ij}\in \mrm{Sym}^2(T^*M)$, consider the variation $g(s)=g+sh$ and let
\begin{equation}
\div_f h=e^{f}\div( e^{-f} h)=\div h-h (\nabla f, \cdot),
\end{equation}
$\div_f^{\dagger}$ be the adjoint of $\div_f$ with respect to
the weighted $L^2$-inner product
\begin{equation}
(\cdot, \cdot)_f=\int_M <\cdot, \cdot> e^{-f}dV,
\end{equation}
\begin{equation}
\Delta_f h:= \Delta h-\nabla f \cdot \nabla h,
\end{equation}
and
\begin{equation} \label{cLf}
\cL_f h=\frac{1}{2}\Delta_f h+Rm(h,\cdot)=\frac{1}{2}\Delta_f h_{ik} +R_{ijkl}h_{jl}.
\end{equation}
Then the second variation $\delta^2_g\nu(h,h)$ of the $\nu$-entropy  is given in \cite{CZ12} by
\begin{align*}
\delta^2_g\nu(h,h)=\left.\frac{d^2}{ds^2}\right|_{s=0}\nu(g(s))
=\frac{1}{(4\pi\tau)^{n/2}}\int_M  <N_f h, h> e^{-f} dV,\\
\end{align*} where the {\it Jacobi operator} (also known as the {\it stability operator}) $N_f$ is defined by
\begin{equation}
N_fh:=\cL_f h +\div_f^{\dagger}\div_f
h+ \frac{1}{2}\nabla^2\hat v_h-Rc \ \frac{\int_M <Rc,
h>e^{-f}\,dV}{\int_M Re^{-f}\,dV},
\end{equation}
and $\hat v_h$ is the unique solution of
\begin{align*}
\Delta_f \hat v_h+\frac{\hat v_h}{2\tau}=\div_f\div_f
h,\qquad\int_M \hat v_h e^{-f}\, dV=0.
\end{align*}
For more details, we refer the reader to our previous paper \cite{CZ12} or Section 2 below.
Note that ${Sym}^2(T^*M)$ admits the following standard direct sum decomposition:
\be\label{decomposition0}
\mrm{Sym}^2(T^*M) = \mrm{Im}(\div_f^\dagger) \oplus\Ker (\div_f).
\ee
The first factor
\begin{equation*}
\mrm{Im}(\div_f^\dagger)= \{ \div_f^\dagger (\omega) \ | \ \omega \in \Omega^1(M)\}\\
= \{\mathscr{L}_X g \ |\ X=\omega^{\sharp} \in \mathscr{X} (M)\}
\end{equation*}
represents deformations $g(s)$ of $g$ by diffeomorphisms. Since the $\nu$-entropy is invariant under diffeomorphisms, the second variation vanishes on this factor.

In \cite{CZ12}, we observed that $\div_f(Rc)=0$ and showed that % $N_f=\cL_f$ on $\Ker(\div_f)_{0}$, and
$Rc$ is an eigen-tensor of $\cL_f$ with eigenvalue\footnote{Note the different sign convention we used in \cite{CZ12} for eigenvalues of $\cL_f$: In \cite{CZ12}, $\lambda$ is an eigenvalue of $\cL_f$ if $-\cL_f h=\lambda h$ for some symmetric 2-tensor $h\neq 0$.} $1/2\tau$, i.e., $\cL_f Rc=\frac 1{2\tau} Rc$. Moreover, for any linearly stable compact shrinking Ricci soliton, we proved that $1/2\tau$ is the only positive eigenvalue of $\cL_f$ on $\Ker(\div_f)$ with multiplicity one. % Here an eigenvalue of $L_f$ means a real number $\lambda$ such that $L_f h = \lambda h$ for some nonvanishing symmetric two tensor $h$.
Very recently, Mehrmohamadi and Razavi \cite{MR} made some new progress. In particular, they showed that $N_f$ vanishes on $\mrm{Im}(\div_f^\dagger)$, extending a similar result in \cite{CHI04, CH15} for positive Einstein manifolds to the compact shrinking Ricci soliton case. In addition, in terms of the operator $\cL_f$, they showed that (i) if a compact shrinking Ricci soliton $(M^n, g, f)$ is linearly stable, then the eigenvalues of $\cL_f$ on $\mrm{Sym}^2(T^*M)$, other than $\frac 1 {2\tau}$ with multiplicity one,
must be less than or equal to $\frac{1}{4\tau}$; (ii) if a compact shrinking soliton $(M^n, g, f)$ has $\cL_f\leq 0$ on $\mrm{Sym}^2(T^*M)$, except on scalar multiples of $Rc$, then $(M^n, g, f)$ is linearly stable (see Theorem 1.3 and Theorem 1.4 in \cite{MR}, respectively).

Clearly, the nonpositivity of the second variation of $\nu$, i.e., $\delta^2_g\nu(h,h)\le 0$, is implied by the nonpositivity of the stability operator $N_f$ on the space $\mrm{Sym}^2(T^*M)$ of symmetric 2-tensors. Thus, studying linear stability of compact shrinking Ricci solitons requires a closer look into the eigenvalues and eigenspaces of $N_f$, especially its leading term  $\mathcal{L}_f$ defined by \eqref{cLf},  acting on $\mrm{Sym}^2(T^*M)$.
Since $\div_f (Rc)=0$, we can further decompose $\Ker (\div_f)$ as
\[\Ker (\div_f) =\bR\cdot \mrm{Rc} \oplus\Ker (\div_f)_{0},\]
where $\bR\cdot \mrm{Rc} =\{\rho Rc \ | \ \rho\in \bR\}$ is the one dimensional subspace generated by the Ricci tensor $Rc$, and
\be\label{perp}
\Ker(\div_f)_{0}=\{h\in \Ker(\div_f) \ | \ \int_M <h, Rc> e^{-f}\, dV=0\}
\ee
denotes the orthogonal complement of $\bR\cdot \mrm{Rc}$ in $\Ker(\div_f)$ with respect to the weighted inner product (1.3). % $\int_M <\cdot, \cdot> e^{-f}\, dV$.
Accordingly, we can refine the decomposition of $\mrm{Sym}^2(T^*M)$ in \eqref {decomposition0} by
\be\label{decomposition}
\mrm{Sym}^2(T^*M) = \mrm{Im} (\div_f^\dagger) \oplus \bR\cdot \mrm{Rc} \oplus\Ker (\div_f)_{0}.
\ee

\smallskip
In this paper,  by exploring decomposition \eqref{decomposition}, we are able to further improve our previous work in \cite{CZ12} and the work of Mehrmohamadi and Razavi \cite{MR}. Our main results are as follows.

\begin{theorem}\label{main theorem} Let  $(M^n, g, f)$ be a compact shrinking Ricci soliton satisfying equation (1.1). Then,

\begin {itemize}

\smallskip
\item[(i)]  the decomposition of $\mrm{Sym}^2(T^*M) $ in \eqref{decomposition}
is both invariant under $\cL_f$  and orthogonal with respect to the second variation $\delta^2_g\nu$ of the $\nu$-entropy.

\smallskip
\item[(ii)] the eigenvalues of $\cL_f$ on $\mrm{Im}(\div_f^\dagger)$ are strictly less than $\frac 1{4\tau}$.

\end {itemize}
\end{theorem}

\begin{theorem}\label{main theorem2}
A compact shrinking Ricci soliton $(M^n, g, f)$ is linearly stable if and only if $\cL_f\leq 0$ on $\Ker(\div_f)_{0}$.
\end{theorem}

\begin{remark}\label{rmk1.2}
Theorem \ref{main theorem} and Theorem \ref{main theorem2} above are extensions of similar results by Hamilton, Ilmanen and the first author in \cite{CHI04} (see also Theorem 1.1 in \cite{CH15}) for positive Einstein manifolds. % with a constant potential function $f$; see Section 4 for more details.
\end{remark}

While there have been a lot of progress in recent years in understanding geometry of general higher dimensional ($n\ge 4$) complete noncompact gradient shrinking Ricci solitons, especially in dimension four, e.g., \cite {CaoZhou10, CN09, Chen09, CLY11, KW15, MW15, MW17} and \cite {CDS19, BCCD22}, very little is known about the geometry of general compact shrinking Ricci solitons in dimension $n=4$ or higher. On the other hand, for possible applications of the Ricci flow to topology, one is mostly interested in the classification of stable shrinking solitons, since unstable ones could be perturbed away hence may not represent generic singularities of the Ricci flow. Thus, exploring the variational structure of compact Ricci shrinkers becomes rather significant.

We point out that Hall and Murphy \cite{HM} have proven that compact shrinking K\"ahler-Ricci solitons with Hodge number $h^{1,1}>1$ are unstable, thus extending the result of Cao-Hamilton-Ilmanen \cite{CHI04} for Fano K\"ahler-Einstein manifolds to the shrinking K\"ahler-Ricci soliton case.
In particular, the Cao-Koiso soliton on $\CP^2\#(-\CP^2)$ and Wang-Zhu soliton on $\CP^2\#(-2\CP^2)$ are unstable.  In addition, Hall-Haslhofer-Siepmann \cite{HHS} and Hall-Murphy \cite{HM14} have shown that the Page metric \cite{Page} on $\CP^2\#(-\CP^2)$ is unstable.
Most recently,  Biquard and Ozuch \cite{BO} proved that the Chen-LeBrun-Weber metric  \cite{CLW08} on $\CP^2\#(-2\CP^2)$ is also unstable. We hope our new results in this paper will play a significant role in future study of linear stability of shrinking Ricci solitons, especially in classifying compact $4$-dimensional linearly stable shrinking Ricci solitons.

\bigskip

\noindent {\bf Acknowledgements.} We thank Professor Detang Zhou for his interest in our work and for his helpful question that led to a simpler proof of Proposition  \ref{thm(ii)}.

\medskip
\section{Preliminaries}

In this section, we fix our notation and recall some useful facts that will be used in the proof of Theorem \ref{main theorem}.
First of all, by scaling the metric $g$, we may assume that  $\tau=1$ in equation (1.1) so that
\be\label{soliton}
Rc+\d^2f =\frac{1}{2}g.
\ee
We also normalize $f$ so that
$$(4\pi)^{-\frac{n}{2}}\int_M e^{-f}\,dV=1.$$
From now on, we shall assume that $(M^n, g, f)$ is a compact shrinking Ricci soliton satisfying \eqref{soliton}.

As in \cite{CZ12}, for any symmetric 2-tensor $h=h_{ij}$ and 1-form
$\omega=\omega_i$, we denote
$$\text{div}\:\omega:=\nabla_i\omega_i, \qquad (\text{div}\:
h)_i:=\nabla_jh_{ji}.$$
Moreover, as done in \cite{Cao08b, CZ12}, we
define $ \div_f (\cdot) := e^{f}\div( e^{-f} (\cdot))$, or more specifically,
\begin{equation} \label{eq1.2}
\div_f \omega =
\div\omega -\omega (\nabla f)=\nabla_i\omega_i-\omega_i\nabla_if,
\end{equation}
and
\begin{equation} \label{eq1.3}
\div_f h=\div h-h (\nabla f, \cdot)=\nabla_jh_{ij}-h_{ij}\nabla_jf.
\end{equation}
We also define the operator $\div_f^{\dagger}$ on functions by
\be
\div_f^{\dagger} u=-\d u, \qquad u\in C^{\infty}(M)
\ee
and on 1-forms by
\begin{equation} \label{div*wLie}
(\div_f^{\dagger} \omega)_{ij}=-\frac 1 2 (\nabla_i\omega_j+\nabla_j\omega_i)
=-\frac 12 \mathscr{L}_{\omega^{\sharp}}g_{ij},
\end{equation}
where $\omega^{\sharp}$ is the vector field dual to $\omega$ and $\mathscr{L}$ denotes the Lie derivative,
so that
\begin{equation} \label{eq1.5}
\int_M e^{-f} <\div_f^{\dagger}\omega, h>dV= \int_M e^{-f}<\omega, \div_f h>dV,
\end{equation}
for any symmetric 2-tensor $h$.

 Clearly, $\div_f^{\dagger}$ is just the adjoint of $\div_f$ with respect to
the weighted $L^2$-inner product
\begin{equation}\label{eq1.6}
(\cdot, \cdot)_f=\int_M <\cdot, \cdot> e^{-f}dV.
\end{equation}
\begin{remark} If we denote by $\div^{*}$ the adjoint of $\div$ with respect to the usual $L^2$-inner product
\be\label{no weight}
(\cdot, \cdot)=\int_M <\cdot, \cdot>dV,
\ee
 then, as pointed out in \cite{Cao08b}, one can easily verify that
\begin{equation} \label{eq1.8}
\div_f^{\dagger}=\div^{*}.
\end{equation}
\end{remark}

\noindent Finally, we denote
\begin{equation} \label{eq1.7}
\Delta_f :=e^f\div(e^{-f}\d)=\Delta -\nabla f \cdot \nabla,
\end{equation}
which is self-adjoint with respect to the weighted $L^2$-inner product (\ref {eq1.6}),
$$Rm(h,\cdot)_{ik}:=R_{ijkl}h_{jl},$$
and define the operator
\be
\cL_f h =\frac 1 2\Delta_f h + Rm(h,\cdot)
\ee
on the space of symmetric 2-tensors. It is easy to see that, like $\Delta_f$,  $\mathcal{L}_f $ is a self-adjoint operator
with respect to the weighted $L^2$-inner product  (\ref{eq1.6}).
\smallskip

Now we restate the second variation of the $\nu$-entropy derived in \cite{CZ12} with $\tau=1$.

\begin{theorem}\label{thm2.1} {\bf (\cite{CZ12})}
Let $(M^n, g, f)$ be a compact shrinking Ricci soliton satisfying \eqref{soliton}. For any symmetric
2-tensor $h=h_{ij}$, consider the variation
$g(s)=g_{ij}+sh_{ij}$. Then the second variation
$\delta^2_g\nu(h,h)$ is given by
\be\label{second variation}
\delta^2_g\nu(h,h)=\left.\frac{d^2}{ds^2}\right|_{s=0}\nu(g(s))
=\frac{1}{(4\pi)^{n/2}}\int_M  <N_f h, h> e^{-f} dV,\\
\ee
where the stability operator $N_f$ is given by
\begin{equation}\label{eq1.9}
N_f h:=\cL_f h+\div_f^{\dagger}\div_f
h+ \frac{1}{2}\nabla^2\hat v_h-Rc \ \frac{\int_M <Rc,
h>e^{-f}\, dV}{\int_M Re^{-f}\, dV},
\end{equation}
and the function $\hat v_h$ is the unique solution of
\begin{align} \label {eq1.10}
\Delta_f \hat v_h+\frac{\hat v_h}{2}=\div_f\div_f
h,\qquad\int_M \hat v_h e^{-f}\, dV=0.
\end{align}
\end{theorem}

Next,  we recall the following facts (see, e.g., Lemma 3.1 and Lemma 3.2 in \cite{CZ12}).

\begin{lemma} {\bf (\cite{CZ12})} \label{lem3.3} Let $(M^n, g, f)$ be a
compact  shrinking Ricci soliton satisfying \eqref{soliton}. Then,

\begin {itemize}
\smallskip
\item[(i)] $ Rc \in \Ker(\div_f);$

\smallskip
\item[(ii)] $\mathcal{L}_f (Rc)=\frac 12 Rc.$

\end {itemize}
\end{lemma}

We shall also need the following useful identities found by Mehrmohamadi-Razavi \cite{MR}; see also Colding and Minicozzi \cite{CM2}, in which they derived more general versions of identities \eqref{dDu}-\eqref{divLie} that are valid for smooth metric measure spaces.

\begin{lemma}{\bf {(\cite{MR, CM2})}}\label{lem identities}  Let $(M^n, g, f)$ be a
compact shrinking Ricci soliton satisfying \eqref{soliton}. Then, for any function $u$, 1-form $\omega$ and symmetric 2-tensor $h$, the following identities hold
\be\label{dDu}
\d \Delta_f u = \Delta_f \d u - \frac{1}{2}\d u,
\ee
\be\label{divDw}
\div_f\Delta_f \omega = \Delta_f\div_f \omega + \frac{1}{2}\div_f \omega,
\ee
\be\label{div*Dw}
\div_f^\dagger\Delta_f \omega = 2\cL_f\div_f^\dagger \omega - \frac{1}{2}\div_f^\dagger \omega,
\ee
\be\label{LLie}
2\cL_f(\mathscr{L}_{\omega^{\sharp}}g)=\mathscr{L}_{(\Delta_f \omega)^{\sharp}}g + \frac{1}{2}\mathscr{L}_{\omega^{\sharp}}g,
\ee
\be\label{divLh}
2\div_f \cL_f h=\Delta_f\div_f h +\frac{1}{2}\div_f h,
\ee
\be\label{divLie}
\div_f (\mathscr{L}_{\omega^{\sharp}}g)=-2\div_f\div_f^\dagger \omega=\Delta_f \omega +\d(\div_f \omega)+\frac{1}{2}\omega.
\ee
\end{lemma}

For the readers' convenience and the sake of completeness, we provide a quick proof here.

\begin{proof} The above identities follow from direct computations given below.

\smallskip
$\bullet$ For \eqref{dDu}:
\[
\al
\d_i \Delta_f u=&\d_i\d_j\d_ju-\d_i\d_j f\d_j u-\d_j f\d_i\d_j u\\
=& \Delta\d_i u + R_{ijjk}\d_k u - \frac{1}{2}\d_i u + R_{ij}\d_j u -\d_j f\d_j\d_i u\\
=& \Delta_f \d_i u -\frac{1}{2}\d_i u.\\
\eal
\]

$\bullet$ For \eqref{divDw}: It follows from \eqref{dDu} that
\[
\al
\int_M u \div_f(\Delta_f \omega)\, e^{-f} dV =& \int_M -<\Delta_f \d u, \omega>e^{-f}\, dV\\
=& \int_M -<\d (\Delta_f u) + \frac{1}{2}\d u, \omega>\, e^{-f}dV\\
=& \int_M u(\Delta_f\div_f \omega+\frac{1}{2}\div_f \omega)\, e^{-f}dV.\\
\eal
\]

$\bullet$ For  \eqref{div*Dw}:
\[
2\cL_f\div_f^\dagger \omega = -\frac{1}{2}\Delta_f(\d_i \omega_j+\d_j \omega_i)- R_{ikjl}(\d_k \omega_l + \d_l \omega_k).
\]
Notice that
\[
\al
\Delta_f\d_i \omega_j=& \d_k\d_k\d_i \omega_j - \d_k f\d_k\d_i \omega_j\\
=& \d_k(\d_i\d_k \omega_j+ R_{kijl}\omega_l)-\d_k f (\d_i\d_k \omega_j+ R_{kijl}\omega_l)\\
=& \d_i \Delta \omega_j + R_{il}\d_l \omega_j+R_{kijl}\d_k \omega_l + \d_kR_{kijl}\omega_l+R_{kijl}\d_k\omega_l\\
\ & -\d_i(\d_kf\d_k\omega_j) +\d_i\d_kf\d_k\omega_j+R_{kijl}\d_kf\omega_l\\
=& \d_i \Delta_f \omega_j - 2R_{ikjl}\d_k\omega_l  + \frac{1}{2}\d_i \omega_j.
\eal
\]

$\bullet$ For \eqref{LLie}: According to \eqref{div*wLie}, \eqref{LLie}  is equivalent to \eqref{div*Dw}.

\smallskip
$\bullet$ For \eqref{divLh}: Similar to the proof of \eqref{divDw}, \eqref{divLh} is the adjoint of \eqref{div*Dw} with respect to the inner product \eqref{eq1.6}.

\smallskip
$\bullet$ For  \eqref{divLie}:
\[
\al
\div_f(\mathscr{L}_{\omega^{\sharp}}g)_j=& \d_i(\d_i \omega_j+\d_j \omega_i)-\d_i f(\d_i \omega_j+\d_j \omega_i)\\
=& \Delta_f \omega_j + \d_j \d_i \omega_i + R_{jk} \omega_k - \d_j(\d_i f \omega_i)+\d_j\d_i f\omega_i\\
=& \Delta_f \omega_j + \d_j \div_f \omega + \frac{1}{2}\omega_j.\\
\eal
\]
\end{proof}

\begin{remark} Some of the identities in Lemma \ref{lem identities} were first obtained  in \cite{CH15} for positive Einstein manifolds. % Identities \eqref{dDu} and \eqref{divDw} were first shown in  \cite{CH15} for positve Einstein manifolds. Moreover, more general versions of identities \eqref{dDu}-\eqref{divLie} were obtained in \cite{CM2}.

\smallskip
For positive Einstein manifolds, C. He and the first author also showed in \cite{CH15} that the restriction of $N_f$ to the subspace $\mrm{Im} (\div_f^\dagger)$ is zero, i.e., $\left. N_f\right|_{\mrm{Im} (\div_f^\dagger)}=0,$ a fact first noted in Cao-Hamilton-Ilmanen \cite{CHI04}.  By using identities \eqref{divDw}, \eqref {LLie} and \eqref{divLie} in Lemma \ref {lem identities}, Mehrmohamadi and Razavi \cite{MR} were able to generalize this to the case of compact shrinking Ricci solitons.
\end{remark}

\begin{lemma}{\bf (\cite{MR})}\label{lem Nimdiv*=0} Let $(M^n, g, f)$ be a
compact shrinking Ricci soliton satisfying \eqref{soliton}. Then, we have
\[N_f|_{\mrm{Im} (\div_f^\dagger)}=0.\]
\end{lemma}

\proof Notice that, according to \eqref{divLie} and \eqref{divDw},
\be\label{v_Lie}
\al
\div_f\div_f(\mathscr{L}_{\omega^{\sharp}}g)=&\div_f(\Delta_f \omega+\d\div_f \omega+\frac{1}{2} \omega)\\
=& 2\Delta_f (\div_f \omega) + \div_f \omega.
\eal
\ee
Thus, if we denote by $\xi=\mathscr{L}_{\omega^{\sharp}}g$, then according to (\ref {eq1.10})
\be\label {hat v}
\hat{v}_{\xi}=2\div_f \omega.
\ee

Now,  by (\ref{div*wLie}), \eqref{LLie}, \eqref{divLie} and \eqref{hat v}, we obtain
\be\label{Ndiv*}
\al
-2N_f (\div_f^\dagger \omega)=& N_f (\mathscr{L}_{\omega^{\sharp}}g)\\
=& \cL_f(\mathscr{L}_{\omega^{\sharp}}g) + \div_f^\dagger\div_f(\mathscr{L}_{\omega^{\sharp}}g) + \d^2 (\div_f \omega)\\
= & \frac{1}{2}\mathscr{L}_{(\Delta_f \omega)^\sharp}g + \frac{1}{4}\mathscr{L}_{\omega^{\sharp}}g + \div_f^\dagger(\Delta_f \omega) + \div_f^\dagger \left(\d(\div_f \omega)\right) \\
\ & + \frac{1}{2}\div_f^\dagger \omega + \d^2(\div_f \omega)\\
=&0.
\eal
\ee
\qed

\section{Proof of the Main Theorems}

In this section, we prove Theorem \ref{main theorem} and Theorem \ref{main theorem2} stated in the introduction.
Once again, by scaling the metric $g$, we normalize $\tau=1$ and assume that $(M^n, g, f)$ is a compact shrinking Ricci soliton satisfying
\be \label{soliton1}
Rc+\d^2f =\frac{1}{2}g.
\ee

First of all, recall that we have the following direct sum decomposition % for $\mrm{Sym}^2(T^*M)$,

\begin{align} %\label{decomposition}
\mrm{Sym}^2(T^*M)  % = & \ \mrm{Im}(\div_f^\dagger) \oplus\Ker (\div_f) \\
= & \ \mrm{Im} (\div_f^\dagger) \oplus \bR\cdot \mrm{Rc} \oplus\Ker (\div_f)_{0},
\end{align}
where $\bR\cdot \mrm{Rc}$ is the one dimensional subspace generated by the Ricci tensor $Rc$ and $\Ker(\div_f)_{0}$, as defined in \eqref{perp}, denotes the orthogonal complement of $\bR\cdot \mrm{Rc} $ in $\Ker(\div_f)$ with respect to the weighted inner product $\int_M <\cdot, \cdot> e^{-f}\, dV$.

\smallskip
We divide the proof of Theorem \ref{main theorem} into two propositions.

\begin{proposition}\label{invariant subspace}
The subspaces $\mrm{Im}(\div_f^\dagger)$, $\bR\cdot \mrm{Rc}$, $\Ker(\div_f)_{0}$ are invariant subspaces of the linear operator $\cL_f$. Moreover,  (3.2) is an orthogonal decomposition with respect to the quadratic form $\delta^2_g\nu (h,h)$ of the second variation in Theorem \ref{thm2.1}.
\end{proposition}

\proof  Firstly, by \eqref{div*Dw},
\[
\cL_f (\div_f^\dagger \omega) = \frac 12 \div_f^\dagger(\Delta_f \omega +\frac{1}{2} \omega) \in \mrm{Im} (\div_f^\dagger).
\]
This shows that  $\mrm{Im}(\div_f^\dagger)$ is invariant under $\cL_f$.

Next, from Lemma \ref{lem3.3} (ii), we have
\[\cL_f Rc=\frac{1}{2}Rc.\]
Hence, $\bR\cdot {Rc}$ is an invariant subspace of $\cL_f$.

Finally, for any $h\in \Ker(\div_f)_{0}$, it follows from \eqref{divLh} that
\[
\div_f(\cL_f h)=\frac{1}{2}\left(\Delta_f \div_f h+\frac{1}{2}\div_f h\right)=0.
\]
Moreover, since  $\cL_f Rc=\frac{1}{2}Rc$, it follows that % from \eqref{LRc}
\begin{equation*}
\al
\int_M <\cL_f h, Rc> \,e^{-f}dV =&\int_M <h, \cL_f Rc>\, e^{-f}dV \\
=& \frac 1 2 \int_M <h, Rc>\, e^{-f}dV=0,
\eal
\end{equation*}
i.e., $\cL_fh\in \Ker(\div_f)_{0}$. Therefore,  $\Ker(\div_f)_{0}$ is also  invariant under $\cL_f$.

Furthermore, the invariant subspace property just demonstrated together with the fact that  $\mrm{Im}(\div_f^\dagger)$, $\bR\cdot \mrm{Rc}$, and $\Ker(\div_f)_{0}$ are mutually orthogonal to each other (with respect to the weighted inner product) immediately imply that
the decomposition \eqref{decomposition} of $\mrm{Sym}^2(T^*M)$ is also orthogonal with respect to the second variation $\delta^2_g\nu(h,h)$ of the $\nu$-entropy. % in \eqref{second variation}.
\qed

\begin{proposition}\label{thm(ii)} Let $(M^n, g, f)$ be a compact shrinking Ricci soliton satisfying (3.1). Then, the eigenvalues of $\cL_f$ on $\mrm{Im}(\div_f^\dagger)$ are strictly less than $\frac {1}{4}.$
\end{proposition}

\proof
Suppose that $\lambda$ is an eigenvalue of $\cL_f$ on $\mrm{Im}(\div_f^\dagger)$, and
$$\cL_f(\mathscr{L}_{\omega^{\sharp}}g)=\lambda \mathscr{L}_{\omega^{\sharp}}g$$
for some $\mathscr{L}_{\omega^{\sharp}}g\equiv-2\div_f^\dagger \omega\in \mrm{Im}(\div_f^\dagger)$ with $\mathscr{L}_{\omega^{\sharp}}g\neq 0$. We need to show $\lambda <\frac {1}{4}.$

Since $N_f =0$ on $\mrm{Im}(\div_f^\dagger)$ by Lemma \ref{lem Nimdiv*=0}, from \eqref{Ndiv*}, we have
\be\label{NLie}
\al
0=& N_f(\mathscr{L}_{\omega^{\sharp}}g) \\
=&\cL_f(\mathscr{L}_{\omega^{\sharp}}g) + \div_f^\dagger\div_f\mathscr{L}_{\omega^{\sharp}}g+\d^2\div_f w\\
=& \lambda\mathscr{L}_{\omega^{\sharp}}g+\div_f^\dagger\div_f\mathscr{L}_{\omega^{\sharp}}g+\d^2\div_f \omega\\
=& -2\lambda\div_f^\dagger \omega -2\div_f^\dagger\div_f\div_f^\dagger \omega-\div_f^\dagger\d\div_f \omega\\
=&-\div_f^\dagger(2\lambda \omega+2\div_f\div_f^\dagger \omega+\d\div_f \omega).
\eal
\ee

\medskip
\noindent
{\bf Claim.} The following identity holds,
\be\label{eigen divw}
\Delta_f\div_f \omega = (2\lambda-1)\div_f \omega.
\ee

\smallskip
Indeed, it follows from \eqref{LLie} that
\[
2\cL_f(\mathscr{L}_{\omega^{\sharp}}g) = \mathscr{L}_{(\Delta_f \omega+\frac{1}{2} \omega)^\sharp}g.
\]
From \eqref{v_Lie}, we know that
\[
\hat v_{\mathscr{L}_{\omega^{\sharp}}g}=2\div_f \omega.
\]
Here, for any symmetric 2-tensor $h$, $\hat v_h$ is given by  (\ref{eq1.10}).
Hence,
\begin{equation*}
\al
2\hat v_{\cL_f(\mathscr{L}_{\omega^{\sharp}}g)}= & \hat v_{\mathscr{L}_{(\Delta_f \omega+\frac{1}{2}w)^\sharp}g}\\
= & 2\div_f(\Delta_f \omega+\frac{1}{2} \omega)\\= & 2(\Delta_f \div_f \omega+\div_f \omega),
\eal
\end{equation*}
where, in the last step above, we have used \eqref{divDw}.

Since
$\cL_f(\mathscr{L}_{\omega^{\sharp}}g)=\lambda\mathscr{L}_{\omega^{\sharp}}g$, we get

\begin{equation*}
\al
\Delta_f \div_f \omega+ \div_f \omega= & \hat v_{\cL_f(\mathscr{L}_{\omega^{\sharp}}g)}\\
=& \lambda \hat v_{\mathscr{L}_{\omega^{\sharp}}g}\\
= & 2\lambda \div_f \omega,
\eal
\end{equation*}
i.e.,
\begin{equation*}\label{eigen div_f omega}
\Delta_f\div_f \omega = (2\lambda-1)\div_f \omega.
\end{equation*}
This proves the Claim.

\smallskip
Now, we divide the rest of our argument into two cases.

\medskip

\noindent
{\bf Case 1: } $\div_f \omega$ is not a constant.
\smallskip

In this case, by the Claim, $\div_f \omega$ is an eigenfunction of $\Delta_f$ with eigenvalue $1-2\lambda$. On the other hand,  from \cite{CZ12}, we know that the first eigenvalue of  $\Delta_f$ is greater than $1/2$. Thus, $1-2\lambda>\frac{1}{2}$; hence $\lambda<\frac{1}{4}$.

\medskip

\noindent
{\bf Case 2: } $\div_f \omega$ is a constant.

\smallskip
In this case, we have
\[
\al
\int_M |\div_f \omega|^2 e^{-f} dV & = \int_M <\omega, \div_f^\dagger\div_f \omega>e^{-f} dV\\
&=-\int_M <\omega, \d(\div_f \omega)>e^{-f} dV\\
& =0.
\eal
\]
It follows that $\div_f \omega=0$. So \eqref{NLie} becomes

\[
\div_f^\dagger(\lambda \omega+\div_f\div_f^\dagger\omega)=0.
\]
Multiplying both sides of the above identity by $\div_f^\dagger\omega$ and integrating yields

\[
\int_M \left(\lambda|\div_f^\dagger\omega|^2+|\div_f\div_f^\dagger\omega|^2\right) e^{-f} dV=0.
\]
Since $\div_f^\dagger\omega=-\frac 1 2 \mathscr{L}_{\omega^{\sharp}}g\ne 0$ by assumption, we have $\lambda\leq 0 <1/4$.

Therefore, we have shown that $\lambda<\frac{1}{4}$. This concludes the proof of Proposition \ref{thm(ii)} and  Theorem 1.1.

\qed\\

Finally,  we are ready to prove Theorem \ref{main theorem2}.

\proof \ By Theorem \ref{thm2.1}, a compact shrinking Ricci soliton $(M^n, g, f)$ is linearly stable if and only if
 \[ \delta^2_g\nu(h,h):=\frac{1}{(4\pi)^{n/2}}\int_M  <N_f h, h> e^{-f} dV \le 0\]
 for every $h\in \mrm{Sym}^2(T^*M)  = \ \mrm{Im} (\div_f^\dagger) \oplus \bR\cdot \mrm{Rc} \oplus\Ker (\div_f)_{0}.$

However, by Theorem \ref{main theorem}(i) (i.e., Proposition \ref{invariant subspace}), we have
\begin{equation*}
\al
 \int_M  <N_f h, h> e^{-f} dV =& \int_M  <N_f h_1, h_1> e^{-f} dV +\int_M  <N_f h_2, h_2> e^{-f} dV \\
\ & +\int_M  <N_f h_{0}, h_{0}> e^{-f} dV \\
= & \int_M  <N_f h_2, h_2> e^{-f} dV  +\int_M  <N_f h_{0}, h_{0}> e^{-f} dV,
\eal
\end{equation*}
where
\[ h=h_1+h_2+h_{0}, \quad {\mbox{with}} \ h_1\in \mrm{Im} (\div_f^\dagger), \ h_2\in \bR\cdot \mrm{Rc}, \ h_{0}\in \Ker (\div_f)_{0},\]
and, in the last equality, we have used the fact that $\delta^2_g\nu(h_1,h_1)=0$ for $h_1\in \mrm{Im} (\div_f^\dagger)$ due to the diffeomorphism invariance of the $\nu$-entropy.  % $N_f$ vanishes on $\mrm{Im}(\div_f^\dagger)$.

\smallskip
On the other hand,
since $\div_f Rc=0$ and $\cL_f Rc=\frac{1}{2}Rc $ by Lemma \ref{lem3.3},  we obtain
\begin{equation*}
\al
N_f(Rc)=& \cL_f Rc-\frac{\int_M|Rc|^2 e^{-f}dV}{\int_M R\, e^{-f}dV}Rc\\
= & \cL_f Rc-\frac{1}{2}Rc =  0,
\eal
\end{equation*}
where we have used the fact that
\[ \int_M|Rc|^2 e^{-f}dV=\frac 12 \int_M R\, e^{-f}dV,\] because the scalar curvature $R$ satisfies the well-known equation $\Delta_fR=R-2|Rc|^2$.
Hence, $N_f=0$ on $\bR\cdot \mrm{Rc}$, and it follows that \[\int_M  <N_f h_2, h_2> e^{-f} dV =0.\]

Also,  as $N_f=\cL_f$ on $\Ker(\div_f)_{0}$, we immediately conclude that
\begin{equation*}
\al
\int_M  <N_f h, h> e^{-f} dV =\int_M  <N_f h_{0}, h_{0}> e^{-f} dV \\
=\int_M  <\cL_f h_{0}, h_{0}> e^{-f} dV.
\eal
\end{equation*}
Therefore, $\delta^2_g\nu(h,h)\le 0$ if and only if
\[\int_M  <\cL_f h_{0}, h_{0}> e^{-f} dV \le 0.\]
This finishes the proof of Theorem  \ref{main theorem2}. \qed\\

\begin{remark} In the proof of Theorem \ref{main theorem2}, if we use Lemma \ref{lem Nimdiv*=0}  instead of Theorem \ref{main theorem} (i) then we would get the following more explicit information about the Jacobi operator $N_f$.

\begin{proposition}\label{Nf deecomposition}
\be
N_f = \begin{cases}
  0,  & \text{on}  \ \mrm{Im} (\div_f^\dagger); \\
  0, & \text{on}  \ \bR\cdot \mrm{Rc}; \\
 \cL_f  & \text{on} \ \! \Ker(\div_f)_{0}.
\end{cases}
\ee
In particular, $N_f \le 0$ on $\mrm{Sym}^2(T^*M)$ if and only if $\cL_f\le 0$ on $\Ker(\div_f)_{0}$.
\end{proposition}
\end{remark}

\begin{remark} Suppose $\xi=\mathscr{L}_{\omega^{\sharp}}g$ is an eigen-tensor of $\cL_f$ for some 1-form $\omega$, with
\[ \cL_f \xi=\lambda \xi.\]
 Then one can show that $\div_f^\dagger\div_f \xi$ and $\d^2\div_f \omega$ are also eigen-tensors of $\cL_f$ with the same eigenvalue, i.e.,
\[ \cL_f (\div_f^\dagger\div_f \xi) =\lambda (\div_f^\dagger\div_f \xi), \]
and
\[\cL_f (\d^2\div_f \omega)=\lambda (\d^2\div_f \omega).\]
Indeed,  if $\cL_f(\xi)=\lambda \xi$ then, by using the identity
\be\label{div*divLh}
\div_f^\dagger\div_f(\cL_f h)=\cL_f(\div_f^\dagger\div_f h)
\ee
shown in \cite{MR}, we have
\[
\al
\cL_f(\div_f^\dagger\div_f\xi)=& \div_f^\dagger\div_f(\cL_f\xi)\\
=& \lambda (\div_f^\dagger\div_f\xi).
\eal
\]
On the other hand, by setting $u= \div_f \omega$ and combining \eqref{eigen divw} with \eqref{LLie} and \eqref{dDu}, we get
\[
\al
2\cL_f(\d^2 u)=& \cL_f(\mathscr{L}_{\d u}g)\\
=& \frac{1}{2}\mathscr{L}_{(\Delta_f (du))^{\sharp}} g+\frac {1}{2}\mathscr{L}_{\frac{1}{2} \d u} g\\
=& \frac{1}{2}\mathscr{L}_{ \d(\Delta_f u+ u)}g\\
=& \frac{1}{2}\mathscr{L}_{2\lambda\d u}g \\
 =& 2\lambda\d^2 u.
\eal
\]
\end{remark}

 To conclude our paper, we mention two open problems.

\bigskip
{\bf Conjecture 1 (Hamilton; 2004 \cite{Cao06, Cao08b})} $\Bbb{S}^ 4$ and $\CP^2$ are the only $\nu$-stable four-dimensional positive Einstein manifolds.

\medskip
{\bf Conjecture 2 (Cao; 2006 \cite{Cao06, Cao08b})} A $\nu$-stable compact shrinking Ricci soliton is necessarily  Einstein, at least in dimension four.

\begin{remark} Besides $\Bbb{S}^ 4$ and $\CP^2$, the other known positive Einstein 4-manifolds are
the K\"ahler-Einstein manifolds $\CP^1\times \CP^1$, $\CP^2\#(-k\CP^2)$ ($3\le k\le 8$), and the (non-K\"ahler Einstein but conformally K\"ahler) Page metric \cite{Page} on $\CP^2\#(-\CP^2)$ and Chen-LeBrun-Weber metric \cite{CLW08} on $\CP^2\#(-2\CP^2)$. Note that, for $n>4$, C. He and the first author \cite{CH15} have found a strictly stable  positive Einstein manifold, other than the round sphere $\Bbb{S}^n$, in dimension $8$.
\end{remark}

\end{document}